# Analysis of Vibration and Thermal of a Modeled Circuit Board of Automated External Defibrillator (AED) Medical Device

Saidi O. Olalere [1]

[1] Georgia Southern University/Department of Mechanical Engineering, Statesboro, United States; Email: olaleresaidi@gmail.com, so04669@georgiasouthern.edu

*Abstract*— **This study was performed to analysis the vibration and thermal changes experienced by an Automated External Defibrillator (AED) Medical Device which is exposed to shocks from patient's reactions, vibration from mobile ambulances and air ambulance, and the heat changes due to the battery component on the circuit board. Basically, AED is made of plastic with the external parts containing majorly the display, button, pad socket, speaker and more while the internal part entails the circuit boards which housed the components like resistors, capacitors, inductors, integrated circuits and more. In this research, the AED was modeled with the Ansys 2020 workbench and calibrated based on static and dynamic loading to verify the static displacement with the first set of five frequencies obtained based on the un-prestressed conditions. With modification, using the prestressed analysis, the next set of frequencies obtained gives an improved result with 0.0003% error difference between each frequency. The modeled Circuit board was used to examine the vibration and dynamic analysis for the rigid board. Likewise, the thermal analysis was conducted on the modeled Circuit board with the heat source as the battery and the rate of dissipation of heat around the board and its effect on the circuit components.**

## I. Introduction

The Automated External Defibrillator (AED) is a life safer to assist those experiencing sudden and life-threatening cardiac arrest (arrhythmias) of ventricular fibrillation and pulseless ventricular tachycardia to analyze the heart's rhythm by delivery electrical shock to reverberate the heart. The AED is made of plastic, and it has external and internal parts. The external part contains the Pad Expiration Window, Latch, Status Indicator, Battery Compartment, and Battery, Electrode holder, Color display, Manual Override button, Shock button, Pad/Electrode Socket, speaker, and IR port. The internal parts include Main board, Display board, Speaker, Display, Shock Discharge Capacitor and Beeper Speaker.

The effect of sudden cardiac arrest is the cause of more than 350,000 deaths in the United States and the way to resuscitate a heart is through the AED. On average, the response time for the first responder is between 8 – 12 minutes while for every delay, the survival rate reduces by 10% approximately (The American Red Cross, 2021). The AED is an easy-to-use instrument which doesn't require any special training. It is placed in public places for easy access in case of emergency or sudden cardiac arrest. During Chicago's Heart Start program for a period of two years, among the 22 persons, 18 were in a cardiac arrhythmia which were treated with AEDs. Of these 18, 11 survived. Of these 11 patients, 6 were treated by bystanders with no previous training in AED use (Sherry L. Caffrey, 2002).

The AED circuit board for the analysis will be a material made of Epoxy FR-4 with length 254mm and width 216mm while the thickness is 0.5mm. The components of the circuit board include Capacitor, Microcontroller, Flash memory, Analog Digital converter, Field Programmable Gate Array (FPGA), Processor, Audio controller, Inductor and more. Based on the design specification and components, a Finite Element Analysis model is established. The governing equation for the experiment is.

$$Mx'' + Cx' + Kx = F(t) \qquad (1)$$

where M, C and K are mass, damping and stiffness matrices respectively.

The goal is to analysis the effect of vibration and thermal experience on the AED based on its operation. The model will undergo various static and dynamic testing of the modeled circuit board. The dynamic and vibration properties will be analyzed for the modeled circuit board for rigid board based on the number and position of the support.

This research will assist in obtaining reliable results when the Defibrillator is used in mobile ambulances experiencing vibration and road bumping. The Automated External Defibrillator has been a lifesaving device which precedes after a Cardiopulmonary Resuscitation (CPR) is performed. The results from AED are important to understand when CPR should be continued or stopped.

The methodology used for the analysis is Finite Element Analysis (FEA) by designing the model using the Ansys Workbench for the circuit board which will contain integrated circuits. The static and dynamic test will be conducted for the model to determine the bending and damping effect respectively.

## II. LITERATURE REVIEW

Deformation experienced in electric components are classified as vibration, shock, and thermal failure. The AED circuit board is a device that experiences failure which can impact the result obtained in resuscitating an individual experiencing cardiac arrest. Generally, most electronic circuit boards experience random vibration instead of ordinary vibration due to external factors within the vibration environment. Most research on electronics is based on high-cycle fatigue to predict fatigue life of component experiencing sinusoidal vibration. Fatigue failure under sinusoidal vibration loading for component by comparing the vibration failure test, FEA, and theoretical test (Y.S.Chen, 2008).

FEA modeling of a PCB's vibration with rigid boundary conditions and comparing the modeling results with the test results using a rigid fixture to identify the PCB dynamic properties only such as natural frequencies (Jingshu Wu, 2002).

For random vibration fatigue, the circuit board research has extended to the soldering by predicting the fatigue life when subjected to random excitation through vibration loading (Pitarresi J.M, 1993). An experimental validated vibration fatigue damage model of a plastic ball grid array solder joint assembly was developed by (Mei-Ling, 2009) to calculate strain and solder joint survival using three-band technique.

## III. MATERIALS AND METHODS

**Model preparation**

The modeled circuit board components used are Capacitors, Microcontroller, Flash Memory, Analog-Digital Converter, FPGA, Processor, Audio Controller, Battery, and Transistor.

The Automated External Defibrillator is a model with different components for analysis. The base plain which is the board has a measurement of 254mm by 216mm.

The capacitors are in cylinder form with length 40mm and diameter of 35mm. The microcontroller is 10mm by 10mm by 1.4mm. The battery designed is 34.5mm in length while its diameter is 17mm.

The components material with the Young Modulus and Poisson ratio were listed below.

| S/N | Component | Materials | Young Modulus | Poisson ratio | Thermal Conductivity |
|---|---|---|---|---|---|
| 1 | Board | FR4 epoxy | 24 GPa | 0.118 | 0.81 W/mK |
| 2 | Capacitor | Tantalum | 175 GPa | 0.34 | 54.4 W/mK |
| 3 | Microcontroller | Copper | 117 GPa | 0.34 | 385 W/mK |
| 4 | Flash memory | Polystyrene | 3250 MPa | 0.34 | 0.033 W/mK |
| 5 | Analog Digital Converter | Silicon | 140 GPa | 0.275 | 150 W/mK |
| 6 | FPGA | Silicon | 140 GPa | 0.275 | 150 W/mK |
| 7 | Processor | Silicon | 140 GPa | 0.275 | 150 W/mK |
| 8 | audio controller | Copper | 117 GPa | 0.34 | 385 W/mK |
| 9 | Battery | Lithium | 31715.884 Pa | 0.355 | 5.4 W/mK |
| 10 | Transistor | Silicon | 140 GPa | 0.275 | 150 W/mK |

Table 1. Numerical Constant of Component

The modeled AED was developed using Ansys 2020 through the workbench. The model circuit board was designed based on specific dimensions of the board and its components. Also, there are four fixed support for the AED which will be fixed to the plastic casing of the Automated External Defibrillator. The FEA is used for the deformation analysis and the thermal effect on the circuit board and its components. The FEA model is presented in Fig 1. The boundary condition was set at the four edges of the modeled circuit board which are fixed as rigid bodies as seen in a typical Automated External Defibrillator.

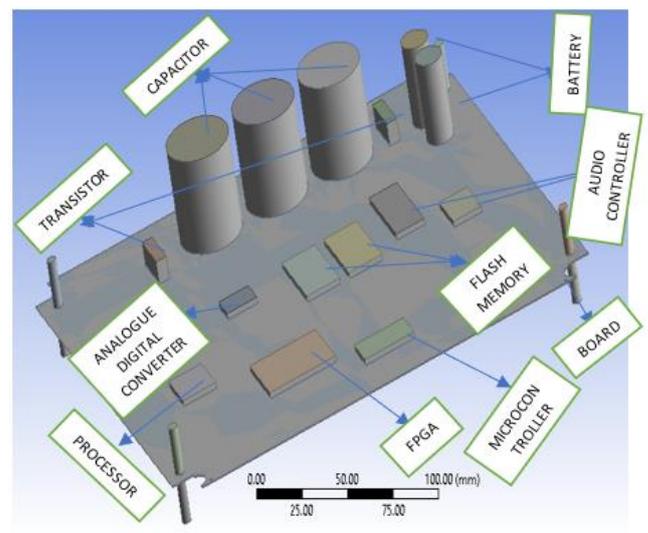

Figure 1. Model of the Circuit Board

**Mesh Selection**

The modeled circuit board was meshed to verify the stress discontinuity of the member components attached to the board. The mesh model produces 90371 nodes and 59671 elements from program-controlled order.

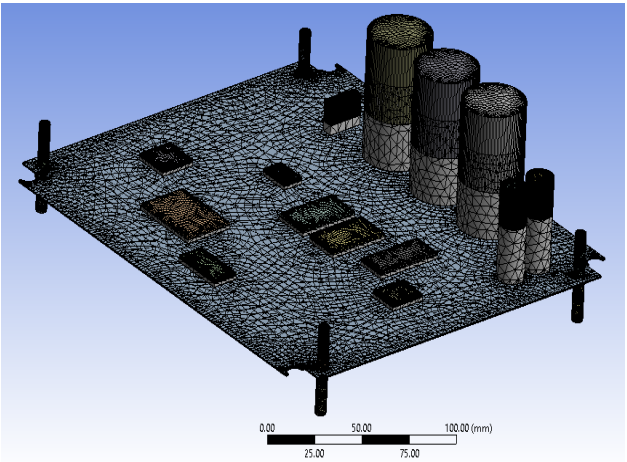

Figure 2. The Mesh of the Circuit Board

## IV. RESULT AND DISCUSSION

**Deformation Analysis for 4 Member Support**

The circuit board undergoes different deformation at different parts of the board. The circuit board experiences the maximum deformation which is more at the middle of the board at a peak of 33.141mm. Likewise, the circuit board experience edge bent along the entire board as shown in Fig 3. The components with high height deform faster which will lead to damage of the circuit board.

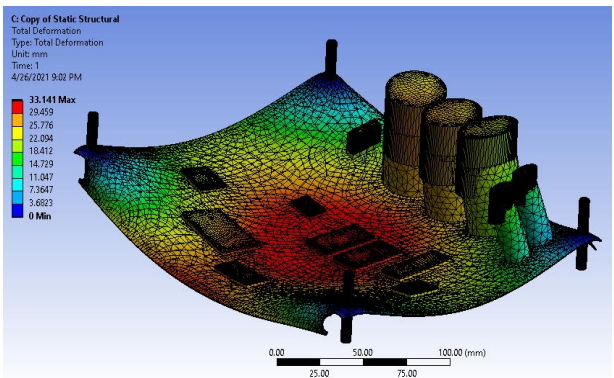

Figure 3. Static Structural Deformation for 4-member Support

The rotational deformation at the Z axis is the largest which transverse is at 10% rate before converging at 60% while at X axis, deformation is experienced at 10 – 60% time series. The Y axis deformation was between 10 – 80% of the time series before finally converging as since in Fig 3.

| Mode | Rotation X | Rotation Y | Rotation Z |
|------|------------|------------|------------|
| 1 | 33.669 | 15.71 | 21.22 |
| 2 | 3.20E+01 | 1.21E+01 | 4.72E+01 |
| 3 | 3.18E+01 | 3.71E+01 | 8.91E+01 |
| 4 | 1.61E-01 | 7.13E+00 | 1.33E-01 |
| 5 | 4.98E+00 | 9.73E-02 | 7.10E-03 |
| 6 | 2.18E-03 | 3.70E-03 | 4.09E-03 |
|  | 102.6497 | 72.17238 | 157.6688 |

Table 2. Mass Participation for 4-member

From Table 2 and Fig 4, slightly above 30% of the effective mass contributed to the mode in the X direction while 50% and 55% effective masses for Y and Z direction respectively.

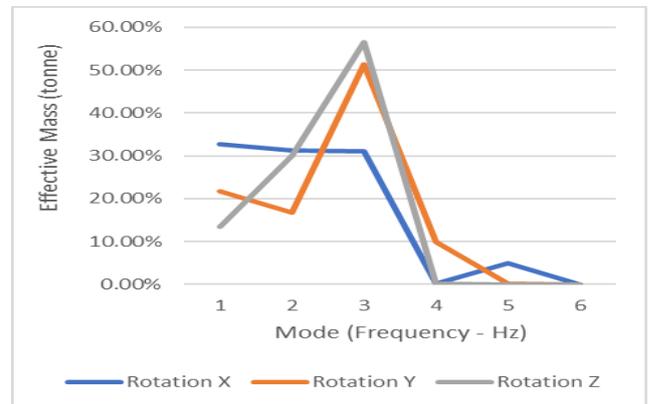

Figure 4. Effective mass to Frequency for 4-member Support

The fatigue experienced by the structure shows that more deformation was visible at the center of the board which has less support compared to the four edges of the modeled circuit board.

In Fig 5, It shows the rate of the deformation at an increasing level which shows the time to failure for the components during the vibration test.

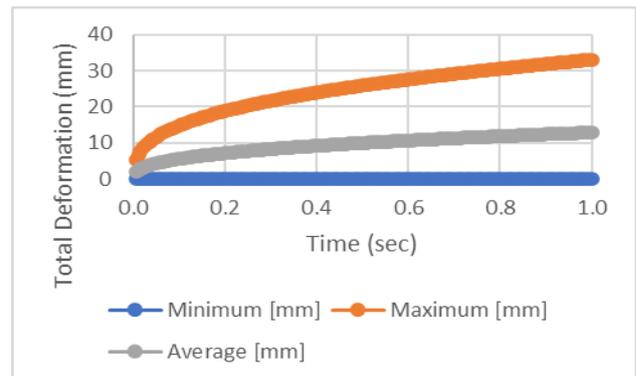

Figure 5. The total deformation rate of the Circuit Board

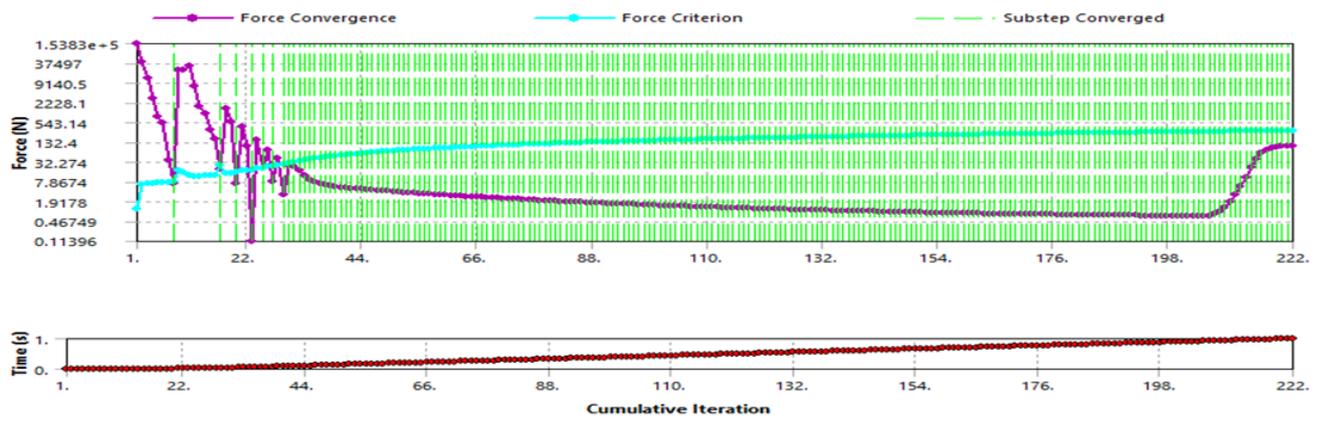

Figure 6. Force Convergence

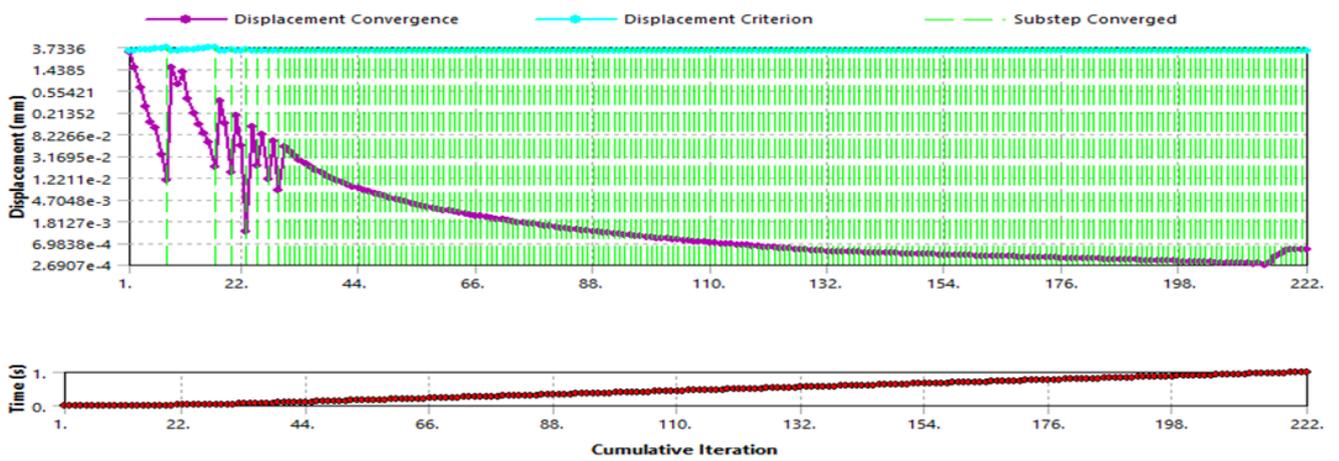

Figure 7. Displacement Divergence

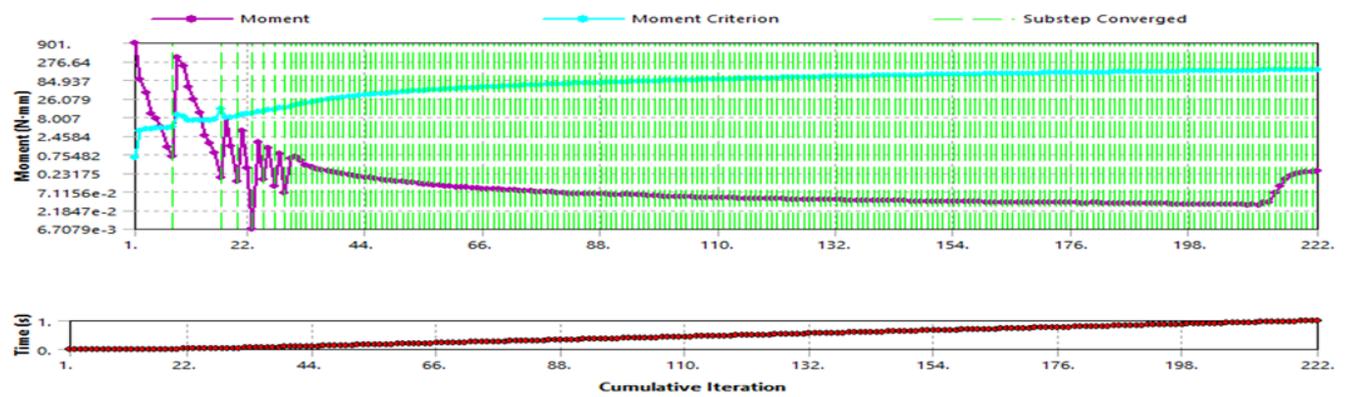

Figure 8. Moment Convergence

For the force convergence in Fig 6, the sub-step converged towards the iteration end point. The convergence experience longer iteration so the load is evenly distributed as seen by the sub-step. At 20% of the loading, normal stiffness experienced is lowered to improve the analysis results.

From the displacement convergence in Fig 7, the normal stiffness was maintained towards the tail end of the analysis signifying that the deformation rate was uniform from 15% iteration.

From the moment convergence in Fig 8, with the mesh refinement which lead to nodes increment shows that the analysis converged with uniformity.

Fig 9 shows the deformation for an un-prestressed condition. The deformation was experienced more at the location of the capacitors which leads to distortion in the circuit board shape.

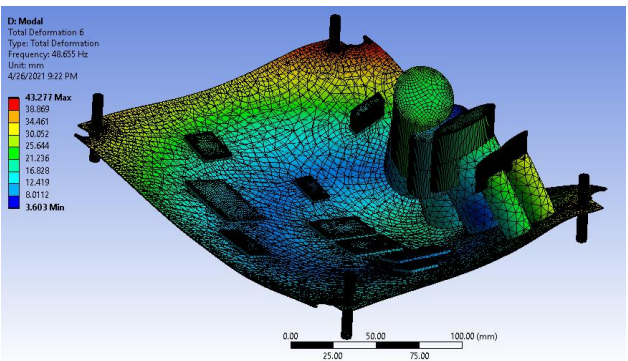

Figure 9. Un-Prestresses Modal Total Deformation

The modal analysis which is used to investigate the vibration on the circuit board is used to evaluate the natural frequencies as shown in Fig 10.

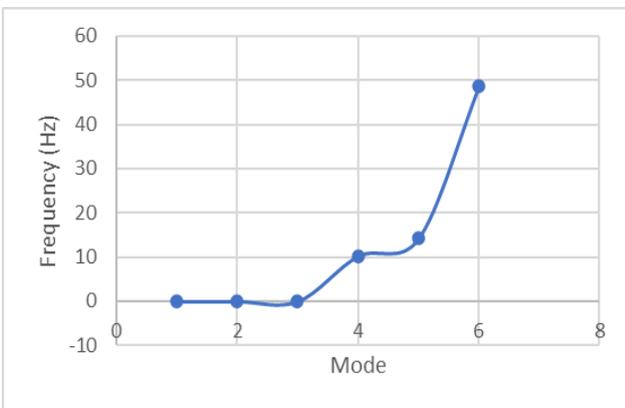

Figure 10. Un-Prestressed Modal Analysis

The prestressed analysis assists in improving the results obtained in the un-prestressed process by modifying the stiffness to reduce the natural frequencies inadequacies. This gives better and improve results for the simulation as the analysis was refined to give better frequencies as against the un-prestressed as seen in Fig 11.

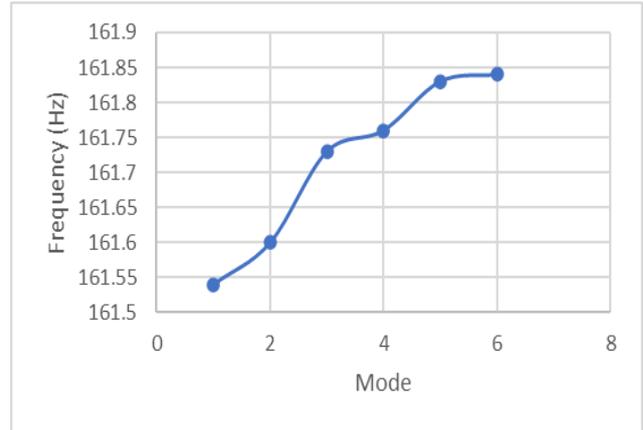

Figure 11. Prestressed Modal Analysis

The modeled circuit board was reinforced with more support members to improve its deformation effect. The deformation peaked at 33.172mm with minimum deformation at the edge of the board as shown in Fig 12.

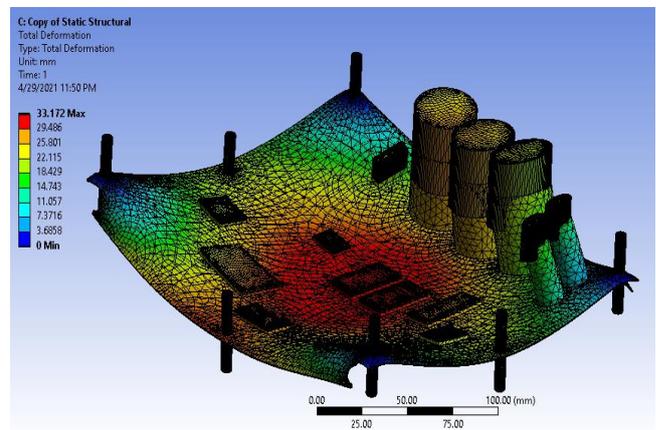

Figure 12. Static Structural Deformation for 8-member Support

| Mode | Rotation X | Rotation Y | Rotation Z |
|---|---|---|---|
| 1 | 33.442 | 25.723 | 5.2324 |
| 2 | 2.67E+00 | 4.11E-01 | 1.53E+02 |
| 3 | 6.17E+01 | 3.95E+01 | 4.72E-01 |
| 4 | 6.09E-03 | 7.47E+00 | 2.01E-01 |
| 5 | 5.24E+00 | 1.06E-02 | 1.01E-03 |
| 6 | 3.00E-01 | 3.94E-03 | 8.40E-03 |
|  | 103.3538 | 73.10313 | 158.784 |

Table 3. Mass Participation for 8-member

From Fig 12 and Table 3, over 95% effective mass was the participating mode in the Z direction, 60% in the X direction and slightly below 60% of effective masses of the mode in the Y direction.

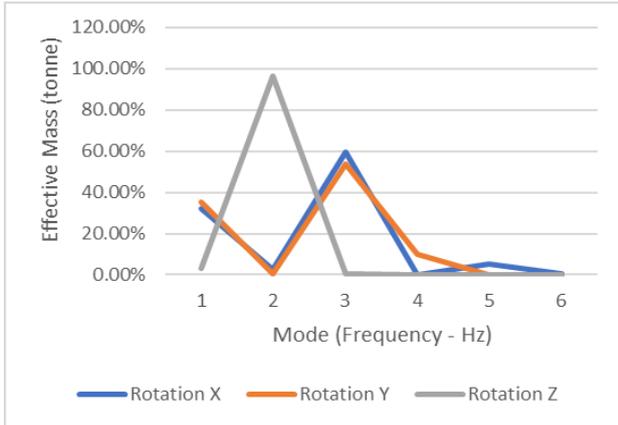

Figure 13. Effective mass to Frequency for 8-member Support

The blue part of the circuit board shown in Fig 14 is the battery of the modeled circuit board which is used to power the board. The temperature is distributed at this point which leads to heat transfer to the entire circuit board. The heat dissipation to the circuit was high, over 40000°c of the initial temperature of 23°c based on the surface area and dissipation rate.

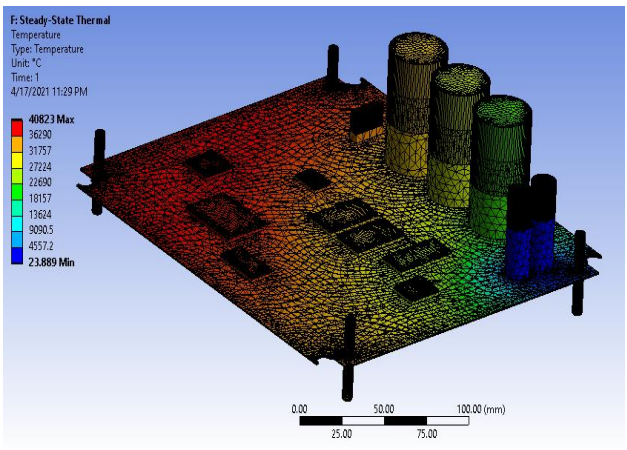

Figure 14. Thermal State of the Circuit Board

The thermal error of the battery as shown in Fig 15 gives a higher value which is due to the temperature dissipated within the components encountered by the circuit board. The temperature effect gives higher thermal error at a small time-interval which will lead to fast rate of the deformation of the circuit board.

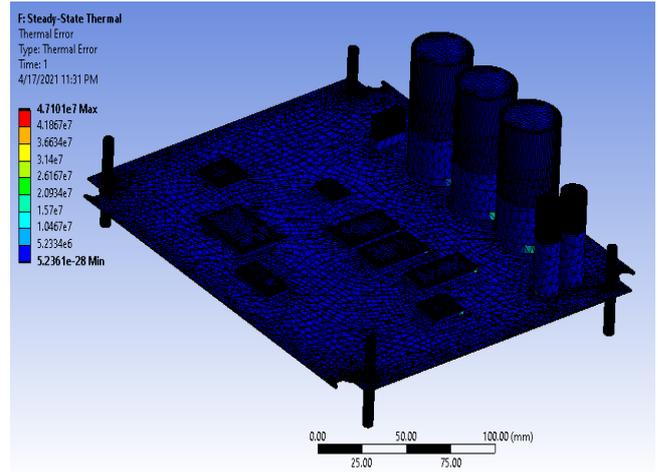

Figure 15. Thermal Error of the Circuit Board

**Discussions**

The analysis was used to carry out the effect of modal analysis and thermal changes on the medical device called Automated External Defibrillator. The experiment was carried out using the Finite Element Analysis in which the model was developed with Ansys workbench to undergo static structure, modal analysis, and steady state thermal.

From other research results, it can be verified that the natural frequencies ranges for both the FEA model and experimental model which may be due to smeared property approach and boundary conditions.

During this research, the detailed joint of the components and the circuit board is not considered but emphasis is based on the face-to-face contact of the components to the circuit board.

Also, the parametric iteration method used to determine the damping varies from 0.001 to 0.005 in a step of 0.005sec.

V. CONCLUSION

The zero frequencies experience is due to the rigid bodies modes which are superfluous, and the weak spring was used to reduce the superfluous effect by minimizing the zero frequencies. The natural frequencies improved when prestressed analysis was used, better results were obtained.

The stiffness of the modeled circuit board is not uniformly distributed because it depends on the components attached to the board and component materials used.

The study shows an improved natural frequency with percentage difference error in 0.0003% shown from the prestressed analysis. The deformation was pronounced more on the component with high heights like capacitors. So, it is suggested that flat capacitors of less height will be suitable in design. Likewise, the heat dissipation by the battery is huge and a better dissipation path is required. Also, lithium battery has a high specific heat capacity of 3582J/(Kg.K) which has high temperature effect on the circuit board. It is suggested that the battery can have a different board or cooling fan should be inbuilt into the Automated External Defibrillator equipment.

ACKNOWLEDGMENT

The author wishes to thank Professor Rahman Mosfequr.